\documentclass[preprint,12pt]{elsarticle}

\usepackage{amssymb}
\usepackage{graphicx}
\usepackage{caption}
\usepackage{subcaption}
\usepackage{xcolor}

\usepackage{amsmath} 
\usepackage{titlesec}
\titleformat{\subsection}[block]{\bfseries}{\thesubsection}{1em}{} 
\titleformat{\subsubsection}[runin]{\bfseries}{\thesubsubsection}{1em}{}
\usepackage{placeins} 
\usepackage{amsmath}
\usepackage{amssymb}
\usepackage{amsthm}
\usepackage{tikz}
\newtheorem{theorem}{Theorem}  
\newtheorem{lemma}{Lemma}
\newtheorem{definition}{Definition}

\newtheorem{corollary}{Corollary}

\newtheorem{conjecture}{Conjecture}

\usepackage{float}
\usepackage{booktabs}
\usepackage{makecell}

\journal{arXiv}

\begin{document}

\begin{frontmatter}



\title{A Survey on Multiset Dimension and Its Variations \\[0.5em] \small 24 April 2026} 




\address[label1]{School of Computer and Information Sciences,
The University of Newcastle, Australia}
\address[label2]{Departments of Mathematics,
University of Indonesia,
Depok - Indonesia}
\address[label3]{Department of Mathematics, Faculty of Science, University of Hafr Al Batin, Saudi Arabia}


\author[label1,label3]{Azzah Albejani}
\ead{azzahahmado.albejani@uon.edu.au}
\author[label1]{Yuqing Lin}
\ead{yuqing.lin@newcastle.edu.au}
\author[label1]{Joe Ryan}
\ead{joe.ryan@newcastle.edu.au}
\author[label2]{Kiki A. Sugeng}
\ead{kiki@sci.ui.ac.id}

\begin{abstract}
The classical notion of metric dimension has led to a wide range of extensions, such as the local, strong, fractional, and 
k-metric dimensions. This naturally raises the question of whether analogous variants can be formulated and studied within the multiset framework. While some progress has been made, particularly on the local multiset dimension, outer multiset dimension, local outer multiset dimension, and edge multiset dimension, however, the area remains far from fully explored. In this paper, we survey the existing variants and consolidate the results currently available in the literature. Furthermore, we identify several directions for future work.

\end{abstract}




\end{frontmatter}




\section{Introduction}
\label{sec:sample1}

The concept of metric dimension in graphs has roots back to 1953 in general metric spaces, though it was formally introduced for graphs by Slater \cite{slater1975leaves}, Harary \& Melter \cite{harary1976metric} independently in the 1970s. Earlier mentions of related ideas also appeared in works by Erd\H{o}s and others \cite{erdos1963two}. 

In addition to its theoretical interest, the work on the metric dimension also has practical significance. For example, Chartrand et al. \cite{chartrand2000resolvability} applied the concept to a real-world problem in the representation and classification of chemical compounds.

The concept of multiset dimension was first introduced by Simanjuntak et al. in 2018 as part of their work to generalization of the metric dimension, the idea is to relax the constraints of metric representations in graphs. Motivation arose from applications in networks where nodes can sense or measure distances to certain reference points but cannot distinguish which specific reference produced which measurement. Since then, the multiset dimension has gained attention as it provides a more flexible framework for node identification in graphs, particularly in contexts such as anonymous networks, robotics, and communication systems.

There are numerous extensions of the metric dimension concept, for example, the local, strong, fractional, and k-metric dimensions, yet corresponding extensions within the multiset-dimension framework remain comparatively limited. In this paper, we summarize the existing results on the multiset dimension and its known variants.

\section{Terms and Definition}

Let $G$ be a connected graph with the vertex set $V(G)$ and edge set $E(G)$, 
the \textit{distance} \( d(u, v) \) between two vertices \( u, v \in V(G) \) is the number of edges in a shortest path between them. The \textit{diameter} of a graph $G$ is the longest distance between any two vertices in $G$. The \textit{edge distance} $d_E(e,f)$ between two edges $e,f \in E(G)$ is the length of a shortest path between vertices $e$ and $f$ in the line graph  $L(G)$.

A $v_0-v_l$ {\it walk} of a graph $G$ is a finite alternating sequence
$v_0, e_1, v_1, e_2, ..., e_l, v_l$
of vertices and edges in $G$ such that $e_i = v_{i-1}v_i$ for each
$i$, $1 \leq i \leq l$. Such a walk may also be denoted by $v_0v_1...v_l$.
We note that there may be repetition of vertices and edges in
a walk. The {\it length} of a walk is the number of edges in the walk.
A $v_0-v_l$ walk is {\it closed} if $v_0=v_l$.
If all the vertices of a $v_0-v_l$ walk are distinct, then the walk
is called a {\it path}. A {\it cycle} is a closed path. A path and a cycle of order \( n \) are denoted by \( P_n \) and \( C_n \), respectively.

A \textit{vertex cover} of a graph is a set $S$ of vertices of $G$ such that every edge of $G$ is incident with at least one vertex of $S$. 
The \textit{vertex covering number} of $G$, denoted by $\alpha(G)$, is the smallest cardinality of a vertex cover of $G$.

A graph is \textit{acyclic} if it has no cycles. A \textit{tree} \( T_n \) is a connected acyclic graph. Any graph without cycles is called a \textit{forest}; thus, the components of a forest are trees. A \textit{star} \( S_n \) graph is a tree consisting of one vertex adjacent to all other vertices.

\textit{Sunflower graph} is a graph obtained by taking a wheel with the central vertex $c$ 
and the $n$-cycle $v_0, v_1, \dots, v_{n-1}$ combined with additional vertices 
$w_0, w_1, \dots, w_{n-1}$, where $w_i$ is joined by edges to $v_i$, $v_{i+1}$, 
where $i + 1$ is taken from modulo $n$.

A \textit{$t$-fold wheel graph} $W_{t,n}$ is a graph derived from a wheel by duplicating the hub vertex one or more times, resulting in $t$ hub vertices, each adjacent to all rim vertices, and not adjacent to each other. 

A \textit{helm graph} $H_n$ is constructed from a wheel $W_n$ by adding $n$ vertices of degree 1, 
one adjacent to each terminal vertex.

A \textit{friendship graph} $f_n$ is a graph consisting of $n$ triangles sharing a common vertex. 

The \textit{$n$-sunlet graph} $S_n$ is obtained by attaching $n$ pendant edges to the cycle $C_n$, having order and size $2n$

The \textit{kayak paddle graph}, denoted by $KP(\vartheta, \lambda, \mu)$, is obtained from two cycles of length $\vartheta \ge 3$ and $\lambda \ge 3$ by joining one vertex of one cycle with a vertex of degree one in a path of length $\mu \ge 2$, and another vertex of the other cycle with the other vertex of degree one in the path of length $\mu$. Its vertex set is 
$
V(KP(\vartheta, \lambda, \mu)) = 
\{\alpha_1, \alpha_2, \ldots, \alpha_{\vartheta}\} 
\cup 
\{\beta_1, \beta_2, \ldots, \beta_{\lambda}\} 
\cup 
\{\gamma_1, \ldots, \gamma_{\mu-1}\},
$
and its edge set is $E(KP(\vartheta, \lambda, \mu)) = 
\{\alpha_i \alpha_{i+1} : 1 \le i \le \vartheta\} \cup \{\beta_j \beta_{j+1} : 1 \le j \le \lambda\} 
\cup \{\gamma_k \gamma_{k+1} : 1 \le k \le \mu - 2\} \cup \{\alpha_1 \gamma_1,\; \gamma_{\mu-1} \beta_1\}.$
where $\alpha_{\vartheta+1} = \alpha_1$ and $\beta_{\lambda+1} = \beta_1$ as defined in \cite{ikhlaq2023new}.

The \textit{dragon graph} $T_{n,m}$ is obtained by joining a vertex $v_n$ of a cycle graph $C_n$
to a vertex $u_1$ of a path graph $P_m$ with a bridge. 
Its vertex set is 
$V(T_{n,m}) = \{v_i, u_j \mid 1 \le i \le n,\; 1 \le j \le m\}$, 
and its edge set is 
$E(T_{n,m}) = \{v_i v_{i+1} : 1 \le i \le n-1\}
\cup \{u_j u_{j+1}  1 \le j \le m-1\}
\cup \{v_1 v_n, v_n u_1\}$ as defined in \cite{ikhlaq2023new}. 

For graph-theoretic terminology not defined in this paper, we follow \cite{bondy2008graph}.

\section{Metric Dimension}

Given a graph $G$, we say a vertex $w$ {\it resolves} a pair of vertices $ u, v $ if $d(u, w) \neq d(v, w)$. For an ordered set of \( k \) vertices $ W = \{w_1, w_2, \ldots, w_k\} \in G$, the {\it representation} of distances of a vertex \( v \) with respect to \( W \) is the ordered \( k \)-tuple 
\[
r(v \mid W) = (d(v, w_1), d(v, w_2), \ldots, d(v, w_k)).
\]
A set of vertices $W \subseteq V(G)$ is a \textit{resolving set} of $G$  if every two vertices of \( G \) have distinct representations. A resolving set with minimum cardinality is called a \textit{metric basis}, and the number of vertices in a metric basis is called the \textit{metric dimension}, denoted by $dim(G)$.

The metric dimension of graphs has been extensively studied over the past decades, resulting in a rich body of literature covering its structural properties, computational complexity, and numerous variants. A comprehensive account of these developments is beyond the scope of this paper. Interested readers are referred to existing survey articles \cite{kuziak2021metric,tillquist2023getting,mohamed2023comprehensive,kratica2014strong} for detailed summary of classical results, algorithmic aspects, and applications of metric dimension. These surveys provide a thorough overview of known bounds, characterizations for special graph classes, and the diverse extensions that have emerged from the original concept.

\section{Metric Dimension Variations}
Over the years, many variations of the classical metric dimension have been introduced. These variants either impose extra structural constraints on the resolving sets, alter the conditions of resolvability, or even aim to identify different types of elements in the graph.

\subsection{The Partition Dimension}
\textit{The partition dimension} of a graph was first introduced in \cite{chartrand2000partition}. It is a variant of the classical metric dimension in which vertices are distinguished based on their distances to subsets of vertices rather than to individual resolving vertex. Formally, a partition $\Pi = \{P_1, P_2, \dots, P_t\}$ of the vertex set $V(G)$ is called a \textit{resolving partition} if every pair of distinct vertices in the graph has a unique vector of minimum distances to the sets in $\Pi$. The smallest number of sets in such a resolving partition defines \textit{the partition dimension} of the graph, denoted by $pd(G)$.


\begin{table}[htbp]
\centering
\renewcommand{\arraystretch}{1.5}
\setlength{\arrayrulewidth}{0.8pt} 
\begin{tabular}{|c|c|}
\hline
\boldmath{$G$} & \boldmath${pd}(G)$\\ 
\hline

Connected graph, $n \ge 4$, not a path or complete graph & $3 \le pd(G) \le n-1$ \cite{chartrand2000partition}\\ \hline
$n$-cycle, $n \ge 3$ & $3$ 
\cite{chartrand2000partition}\\ \hline
$3$-cube $Q_3$ & $3$
\cite{chartrand2000partition}\\ \hline
Petersen graph & $4$
\cite{chartrand2000partition}\\ \hline
Nontrivial connected graph $H$,  $(H \times K_2)$ & $\le pd(H) + 1$
\cite{chartrand2000partition}\\ \hline
$n$-cube $Q_n$, $n \geq 2$ & 
$\leq n+1$ 
\cite{chartrand2000partition}\\ \hline
Complete bipartite graph $K_{r,r}$ &  $\leq r + 1$
\cite{chartrand2000partition}\\ \hline
Complete bipartite graph $K_{r,s}$, $r \neq s$ &  $\leq \max\{r, s\}$
\cite{chartrand2000partition}\\ \hline
$K_{1,n-1}$, $K_n - e$, $K_1 + (K_1 \bigcup K_{n-1})$, $n \geq 3$ & $n - 1$
\cite{chartrand2000partition}\\ \hline

\end{tabular}
\caption{The partition dimension of graphs.}
\label{tab:pd for graphs}
\end{table}
\FloatBarrier

A set of fundamental results on the partition dimension was developed in \cite{chartrand2000partition}, including bounds for bipartite graphs and general connected graphs in terms of their order and diameter; the reader is referred to \cite{ kuziak2021metric, chartrand2000partition} for details.

\subsection{The Edge Metric Dimension}
\textit{The edge metric dimension} extends the idea of metric identification from vertices to edges. While the classical metric dimension focuses on uniquely identifying the vertices of a graph, the edge metric dimension focuses on distinguishing edges through their distances to some vertices.

Formally, given a connected graph $G$, the distance between a vertex $v \in V(G)$ and an edge $e = xy \in E(G)$ is defined as $d_G(e, v) = \min\{d_G(v, x), d_G(v, y)\}$. A set of vertices $S \subset V(G)$ is called an \textit{edge resolving set} if for every pair of distinct edges $e, f \in E(G)$, there exists a vertex $v \in S$ such that $d_G(e, v) \ne d_G(f, v)$. The smallest such set is called an \textit{edge metric basis}, and its cardinality is \textit{the edge metric dimension}, denoted by $\mathrm{edim}(G)$. This parameter was first introduced in \cite{kelenc2018uniquely}.

We have noticed another definition of the edge metric dimension which was given in \cite{article} as follows: 

\begin{definition} 
 Let $W_E=\{e_1,e_2,...,e_k\}$ be an ordered set of the edges of $G$(vertices of $L(G)$) and $f \in E(G)=V(L(G))$, then representation $r(f|W_E)$ of $f$ with respect to $W_E$ is the $k$-tuple $(d_E(f,e_1), d_E(f,e_2), ..., d_E(f,e_k))$. $W_E$ is called edge resolving set if $r(f|W_E)$ is distinct for all $f \in E(G)$. A edge resolving set of minimum cardinality is called an edge basis for $G$ and this cardinality is the edge metric dimension of $G$, denoted by $\dim_E(G)$.  

\end{definition}

This definition uses the edges to resolve the edges, which translate the problem into determine the metric dimension of the line graph of the original graph, thus, in this paper, we follow the notation $\mathrm{edim}(G)$ as defined in \cite{kelenc2018uniquely}.

A number of important results on the edge metric dimension 
were introduced in \cite{kelenc2018uniquely}. These include values for several graph families such as complete bipartite graphs, trees, grids, wheels, and fan graphs, as well as 
NP-completeness and NP-hardness, an $O(\log m)$ approximation bound, structural bounds, and results for certain graph families such as hypercubes. For more details, the reader is referred to \cite{kuziak2021metric}.

\vspace{10pt}

The following table summarize some of the known results:
\begin{table}[!ht]
\centering
\begin{tabular}{|c|c|}
\hline
\boldmath$G$ & \boldmath$\mathrm{edim}(G)$\\ 
\hline
$P_n$ & $1$ \cite{chartrand2000resolvability}\\ \hline
$C_n$ & $2$ \\ \hline
$K_n$ & $n-1$ \cite{bailey2011base}\\ \hline
$T$ & $\sigma(T)-\mathrm{ex}(T)$ \cite{feng2013metric}\\ \hline
$K_{s,t}$ & 
$\begin{cases}
\left\lfloor \frac{2(s+t-1)}{3} \right\rfloor, & \text{if } s \le t \le 2s \\ 
 t-1, & \text{if } t \ge 2s
\end{cases}$ \cite{caceres2007metric}\\ \hline
$W_{1,n}$ & 
$\begin{cases}
3, & n = 3,4 \\ 
4, & n = 5 \\ 
n - \lceil n/3 \rceil, & n \ge 6
\end{cases}$ \cite{eroh2012metric}\\ \hline
$B_n$ & $2n-1$ \cite{eroh2012metric}\\ \hline
$n$-sunlet graph $S_n$ & $ =
\begin{cases}
2, & \text{if $n$ is even},\\
3, & \text{if $n$ is odd}.
\end{cases}$ \cite{article}\\ \hline
$G_{n,k}$ & =
$\begin{cases}
2, & \text{if $1 \le k \le \lceil n/2 \rceil$},\\[2mm]
3, & \text{if $\lceil n/2 \rceil \le k \le n$}.
\end{cases}$ \cite{article} \\ \hline
$D_n$ the family of prism graph & $ 3, \quad \text{for } n \ge 3.$ \cite{article} \\ \hline
\end{tabular}
\caption{Edge metric dimension of graphs.}
\label{tab:edim for graphs}
\end{table}
\FloatBarrier


\subsection{The Local Metric Dimension}
\textit{The local metric dimension} of a connected graph $G$ is defined based on a \textit{local resolving set}, which is a subset of vertices $S \subset V(G)$ such that for every pair of adjacent vertices $u, v \in V(G)$, there exists a vertex $x \in S$ that distinguishes them by their distances; that is, $d_G(u, x) \neq d_G(v, x)$. The smallest such set $S$ is called a \textit{local metric basis}, and its cardinality is known as the \textit{local metric dimension}, denoted $\dim_{\ell}(G)$, $ldim(G)$ \cite{kuziak2021metric}, $\mathrm{lmd}(G)$
\cite{okamoto2010local} or $\beta_l(G)$ \cite{fancy2021local}. In this paper we use the symbol $ldim(G)$. This concept was first introduced in \cite{okamoto2010local}.

The following table summarize some of the known results:

\begin{table}[!ht]
\centering
\renewcommand{\arraystretch}{1.5}
\setlength{\arrayrulewidth}{0.8pt} 
\begin{tabular}{|c|c|}
\hline
\boldmath$G$ & \boldmath${ldim}(G)$\\ 
\hline

Path graph $P_n$ & $1$ \cite{alfarisi2019local}\\ \hline

Cycle graph $C_n$ & $=\begin{cases}
1, & \text{if $n$ is even},\\[2mm]
2, & \text{if $n$ is odd}.
\end{cases}$ \cite{alfarisi2019local}\\ \hline

$K_n$ & $n-1$ \cite{okamoto2010local}\\ \hline
Bipartite graph & $1$ \cite{okamoto2010local}\\ \hline

Complete $k$-partite graph, $k \ge 2$ & $k - 1$ \cite{okamoto2010local}\\ \hline

Wheel graph $W_n$, $n \ge 5$ & $\left\lceil \frac{n}{4} \right\rceil$ \cite{fancy2021local}\\ \hline

Helm graph $H_n$, $n \ge 5$ & $\left\lceil \frac{n}{4} \right\rceil$ \cite{fancy2021local}\\ \hline

Sunflower graph $SF_n$, $n \ge 7$ & $\left\lceil \frac{n}{4} \right\rceil$ \cite{fancy2021local}\\ \hline

\end{tabular}
\caption{The local metric dimension of graphs.}
\label{tab:localdim for graphs}
\end{table}
\FloatBarrier

More results on the local metric dimension can be found in \cite{okamoto2010local}. These include structural properties, bounds in terms of diameter, clique number, and true twin equivalence classes. Further details are available in \cite{kuziak2021metric}.

\subsection{Other Variants of Metric Dimension }

In addition to the above variations of the metric dimension, several other variants have been introduced to capture different aspects of the identification of the vertices based on distances. These include \textit{the mixed metric dimension} which is an extension of both the classical metric dimension and the edge metric dimension, designed to distinguish all elements of the graph, vertices and edges, using a single set of reference vertices, for details, see \cite{kelenc2017mixed}.

The \textit{adjacency dimension} relies on the adjacency representation of vertices, where this adjacency relation allows the resolving set to uniquely distinguish between pairs of vertices based on their direct connections \cite{jannesari2012metric}\cite{kuziak2021metric}. 

\textit{The strong metric dimension} is an extension of the classical metric dimension, developed to capture the true distance relationships within a graph more effectively. A vertex $w$ is said to strongly resolve two vertices $u$ and $v$ if there is a shortest path between $w$ and $u$ containing $v$ or a shortest path between $w$ and $v$ containing $u$ \cite{sebHo2004metric}.

Another variation is \textit{the strong partition dimension} \cite{yero2013strong}, which is based on partitioning the vertex set of a graph instead of selecting individual vertices. The sets of the partition, under the strong resolving condition, act as elements that uniquely determine the position of each vertex.

The \textit{k-metric dimension} generalizes the standard metric dimension by requiring a higher level of redundancy in vertex distinction; every pair of distinct vertices is distinguished by at least $k$ vertices in the resolving set. This concept can be useful in applications demanding fault tolerance and robustness. Notably, when $k=1$, the $k$-metric dimension coincides with the classical metric dimension. This concept was introduced independently in both \cite{adar2017k} and \cite{estrada2013k}.


\textit{The fractional metric dimension} was first described in \cite{arumugam2012fractional}. It is a modified version of the classical metric dimension that allows assigning fractional weights to vertices rather than selecting them completely (as in 0 or 1 choices). 

The \textit{fractional strong metric dimension} assigns a weight between $0$ and $1$ to each vertex instead of selecting vertices entirely. For every pair of distinct vertices, the total weight of vertices that strongly resolve them must be at least $1$, and the goal is to minimize the sum of all vertex weights \cite{kang2013fractional}.


\textit{The $k$-metric antidimension} is a graph parameter introduced to support privacy protection in network analysis, especially in social networks. Unlike the classical metric dimension, which aims to uniquely identify all nodes based on distances to a resolving set, the $k$-metric antidimension aims to maximize anonymity by making some nodes indistinguishable \cite{trujillo2016k}\cite{kuziak2021metric}.

\section{Multiset Dimension}

\textit{The multiset dimension} is an extension of metric dimension where vertices are uniquely identified not by ordered distance vectors, but by multiset of distances to a chosen set of vertices $S$. For a vertex $u$, its multiset representation with respect to $S$, denoted $m(u|S)$, is the multiset of distances from $u$ to each vertex in $S$. A set $S$ is called a \textit{multiset resolving set} if every pair of distinct vertices in the graph has different multiset representations relative to $S$. The minimum size of such a set is \textit{the multiset dimension} of the graph. This parameter can be infinite in some graphs, especially those containing twins vertices that have identical neighborhoods because their multisets of distances to any set $S$ cannot differ \cite{Rinovia2017multiset}.

Let \( G \) be a simple and connected graph with vertex set \( V(G) \). Suppose that \( W \) is a subset of \( V(G) \) and \( v \) is a vertex of \( G \). The representation multiset of \( v \) with respect to \( W \), \( r_m(v \mid W) \), is defined as a multiset of distances between \( v \) and the vertices in \( W \) together with their multiplicities. If \( r_m(u \mid W) \neq r_m(v \mid W) \) for every pair of distinct vertices \( u \) and \( v \), then \( W \) is called a resolving set of \( G \). A resolving set having minimum cardinality is called a multiset basis. If \( G \) has a multiset basis, then its cardinality is called the multiset dimension of \( G \), denoted by \( md(G) \). If \( G \) does not contain a multiset basis, we say that \( G \) has an infinite multiset dimension and we write \( md(G) = \infty \) \cite{Rinovia2017multiset}. 

\section{Known Results on Multiset Dimension }

Several foundational results on the multiset dimension were first established, particularly those concerning existence. In \cite{Rinovia2017multiset}, the followings are given.

\begin{lemma} \cite{Rinovia2017multiset}
No graph has multiset dimension 2.
\end{lemma}

\begin{theorem}\cite{Rinovia2017multiset}
Let $G$ be a graph other than a path. Then ${md}(G) \ge 3$.     
\end{theorem}

\begin{lemma} \cite{Rinovia2017multiset}
${md}(G) \geq \mathrm{dim}(G)$.
\end{lemma}

For positive integers \( n \) and \( d \), let \( f(n, d) \) denote the least positive integer \( k \) for which
\[
\frac{(k + d - 1)!}{k! \, (d - 1)!}+k \geq n.
\]

The following theorem gives a lower bound.

\begin{theorem} \cite{Rinovia2017multiset}
\label{betterlowbound}
If \( G \) is a graph of order \( n \geq 3 \) and diameter \( d \), then 
\[
{md}(G) \geq f(n, d).
\]
\end{theorem}

Further study from \cite{Novi2021some}, Novi et al. have shown that.

\begin{theorem} \cite{Novi2021some}
The multiset dimension is not monotonic to the order of the graph.
\end{theorem}

Furthermore, a construction of a tree is given which lead to the following results. 

\begin{theorem} \cite{Novi2021some}\label{Novitheorem}
Given any positive integer $m \geq 3$, there exists a simple graph whose multiset dimension equals $m$.
\end{theorem}

Following the Theorem \ref{betterlowbound}, for graphs with multiset dimension 3, Novi et al. showed the following bound and also the graph reaching the bound. 
\begin{theorem} \label{nkmd} \cite{Novi2021some}
If \( G \) is a graph on \( n \) vertices of diameter \( k \) with multiset dimension 3, then 
\[
n \leq \frac{k^3 + 3k^2 + 2k + 12}{6},
\]
and this bound is the best possible.
\end{theorem}

Another structural bound relating the multiset dimension to 
the maximum degree of a connected graph is as follows:

\begin{theorem}\cite{alfarisi2020note}
Let $G$ be a connected graphs and let $d$ be the maximum degree of $G$ and $md(G) = k$, we have for $k \ge 3$ and $d < 3k$.
\end{theorem}

Several conditions for a graph to have infinite multiset dimension were also identified. Specifically,

\begin{theorem} \cite{Rinovia2017multiset}
\label{sufficient condition1}
If \( G \) is a graph of diameter at most 2 other than a path, then 
\[
{md}(G) = \infty.
\]
\label{cycles with at most 5 vertices and others}
\end{theorem}
This means that cycles with at most 5 vertices, complete graphs, stars, the Petersen graph, and strongly regular graphs have multiset dimension infinity.

\begin{lemma} \cite{Rinovia2017multiset}
\label{sufficient condition2}
If \( G \) contains a vertex which is adjacent to (at least) three pendant vertices, then 
\[
{md}(G) = \infty.
\]
\end{lemma}

\begin{figure}[ht]
    \centering
    \includegraphics[width=0.4\textwidth]{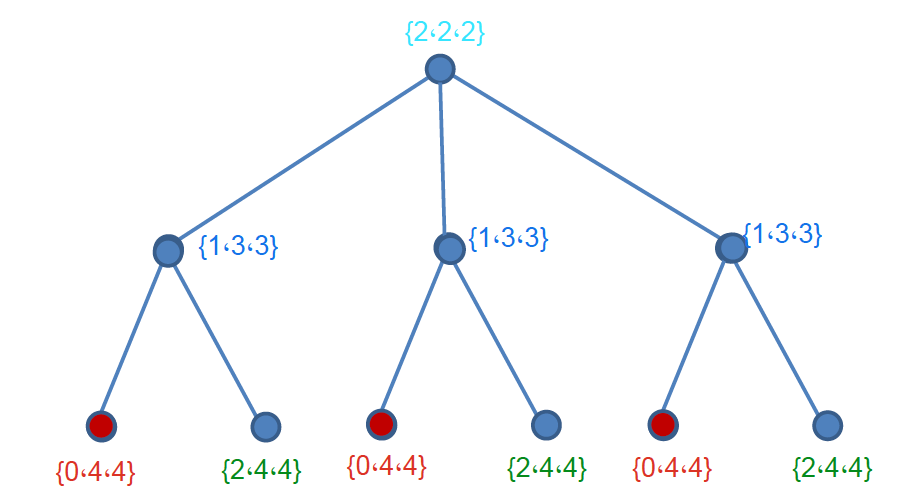}
    \caption{A tree that satisfies neither of the following: diameter at most $2$, nor having a vertex adjacent to at least three pendant vertices, yet still has infinite multiset dimension.}
    \label{fig:sample}
\end{figure}

However, these conditions are neither sufficient nor necessary for a graph to have infinite multiset dimension. For example in Figure~\ref{fig:sample}, there exists a graph that does not satisfy these conditions yet still has \(\operatorname{md}(G) = \infty\).

\FloatBarrier

The following conjecture was proposed. 

\begin{conjecture}\cite{Rinovia2017multiset} \label{conjecture2.1}If $G$ is a graph on $n$ vertices having finite multiset dimension,
then ${md}(G) \le n - 1$.
\end{conjecture}

Conjecture~\ref{conjecture2.1} was proved for trees in \cite{hafidh2019multiset}  

\begin{theorem} \cite{hafidh2019multiset}
Let $T$ be a tree of order $n$ and diameter at least $2$. If $md(T) < \infty$,
then $md(T) \le n - 2$.
\end{theorem}

Furthermore, the following was conjectured.
\begin{conjecture}\cite{hafidh2019multiset}
Let $T$ be a tree. If $md(T) < \infty$, then $md(T) \le n - \operatorname{diam}(T) + 1$
and the bound is sharp.
\end{conjecture}

There are various of results on trees. In \cite{hafidh2019multiset}, the following was proved.  
\begin{theorem} \cite{hafidh2019multiset}
Let $G$ be a lobster. If $P$ is the minimum $2$-center path of $G$, then
the following are equivalent.
\begin{enumerate}
    \item $G$ has finite multiset dimension.
    \item The only component of $G - E(P)$ with infinite multiset dimension is an $S_4$.
    \item If $H$ is a component of $G - E(P)$ then $[H]$ has at most $4$ components which
    are either a $P_2$, a $P_3$, or an $S_4$, with at most two $P_2$s and two $S_4$s.
\end{enumerate}
\end{theorem}

\begin{theorem}\cite{hafidh2019multiset}
Let $G$ be a caterpillar with $P$ its minimum $1$-center-path. $md(G)$ is
finite if and only if every vertex in $P$ has at most $2$ neighbors in $G - P$.
\end{theorem}

In the paper \cite{ikhlaq2023new}, the authors defined the kayak paddle graphs, dragon graphs, and the comb products of two path graphs, and established that all of these graphs have multiset dimension $3$. Marcelo et al.\cite{marcelo2025multiset} have determined the multiset dimension of prisms. Moreover, Marcelo et al.\cite{marcelo2025multiset2} proved that all cylindrical graphs $P_m \square C_n$, 
with $m, n \ge 3$, have finite multiset dimension. For exact values, see Table~\ref{tab:md}.

The concepts of ID-colorings and
ID-graphs were introduced and studied in \cite{chartrand2021distance}. These concepts are equivalent to $m$-resolving sets, where the cardinality of a smallest $m$-resolving set is the multiset dimension. Later, Kono and Zhang studied identification colorings for trees in \cite{kono2021vertex}, for caterpillars in \cite{kono2022note}, and for grids and prisms in \cite{kono2022vertex}. Spectra of grids in \cite{marcelo2025vertex}. Lollipop graphs have been introduced in \cite{cai2025identification}. The paper \cite{wang2024graph} presents some relations between ID-graphs and their ID-indices; give a lower bound on the ID-index of a graph; and determine the ID-indices of paths, grids, cycles, prisms, complete graphs, some complete multipartite graphs, and some caterpillars. Hakanen et al. in \cite{hakanen2024complexity} proved that the multiset resolving sets and the ID-colorings of graphs are equivalent problems. The definition of ID-colorings is given below. 

\begin{definition}\cite{hakanen2024complexity}
Let $G$ be a connected graph of diameter $d$ and a set of vertices $S \subset V(G)$. For every vertex $x \in V(G)$, the \emph{code} of $x$ with respect to $S$ is the $d$-vector 
$\vec d(x|S) = (a_1, a_2, \ldots, a_d)$, where $a_i$, with $i \in \{1, \ldots, d\}$ represents the number of vertices in $S$ at distance $i$ from $x$. If all the codes of vertices of $G$ are pairwise different, then $S$ is called an \emph{identification coloring}
or \emph{ID-coloring}. Moreover, a graph $G$ that has an ID-coloring is called an \emph{ID-graph}. For any ID-graph $G$, the cardinality of a smallest ID-coloring is the \emph{ID-number} of $G$, denoted by $\mathrm{ID}(G)$.
\end{definition}

And the following results is show.
\begin{theorem}\cite{hakanen2024complexity} Let $G$ be a graph of diameter $d$. Then $S \subset V(G)$ is an ID-coloring for $G$ if and only if $S$ is a multiset resolving set for $G$.
\end{theorem}

We have reviewed the literature on ID-coloring, including \cite{chartrand2021distance,kono2021vertex,kono2022note,kono2022vertex,marcelo2025vertex,wang2024graph}. All results are already covered in the multiset dimension setting, except for lollipop graphs \cite{cai2025identification}. The following results establish the multiset metric dimension or derive upper bounds for the corresponding lollipop graphs.

\begin{theorem} \cite{cai2025identification}
The lollipop graph $T_{3,n+1}$ $(n \ge 2)$ has an identification coloring number $r$ for an ID-coloring if and only if $3 \le r \le n+1$.
\end{theorem}

\begin{theorem} \cite{cai2025identification}
The lollipop graph $T_{4,n+1}$ $(n \ge 1)$ has an identification coloring number $r$ for an ID-coloring:
\begin{itemize}
    \item when $n = 1$ or $n = 2$, if and only if $r = 3$;
    \item when $n \ge 3$, if and only if $3 \le r \le n+2$.
\end{itemize}
\end{theorem}

\begin{theorem} \cite{cai2025identification}
The lollipop graph $T_{5,n+1}$ $(n \ge 1)$ has an identification coloring number $r$ for an ID-coloring:
\begin{itemize}
    \item when $n = 1$, if and only if $r = 3$ or $r = 4$;
    \item when $n \ge 2$, if and only if $3 \le r \le n+4$.
\end{itemize}
\end{theorem}


Classifying the problem of computing the multiset dimension as NP-complete \cite{hakanen2024complexity} indicates that finding an optimal solution is computationally very challenging, and that no efficient algorithms are known to solve it in reasonable time for all graphs. As a result, research often focuses on approximation algorithms or on solving the problem for specific graph classes. This motivates our study to examine the multiset dimension in structured graphs such as hypercubes, and to explore whether understanding its relationship with other graph dimensions could lead to more efficient methods or simplifications in special cases.

\noindent
\par
\vspace{5pt}

\begin{table}[!ht]
\centering
\caption{Multiset dimension of some graphs}
\label{tab:md}
\renewcommand{\arraystretch}{1.5}
\setlength{\arrayrulewidth}{0.8pt}
\resizebox{\textwidth}{!}{
\begin{tabular}{|l|c|}
\hline
\textbf{Graph} & \boldmath${md}$ \\
\hline
Path graph $P_n$ & $1$ \cite{Rinovia2017multiset} \\
\hline
Cycle graph $C_n$ & $3$, $n \ge 6$ \cite{Rinovia2017multiset} \\
\hline
Complete graph $K_n$ & $\infty$ \cite{Rinovia2017multiset} \\
\hline
Star graph $S_n$ & $\infty$ \cite{Rinovia2017multiset} \\
\hline
Petersen graph & $\infty$ \cite{Rinovia2017multiset} \\
\hline
Complete $k$-ary tree ($k = 1$ or $2$) & finite \cite{Rinovia2017multiset} \\
\hline
Complete binary tree of height $h$ & $2^h - 1$ \cite{Rinovia2017multiset} \\
\hline
$P_m \square P_n$ & $3$, $m \ge 3$, $n \ge 2$ \cite{Rinovia2017multiset} \\
\hline
$P_n \boxtimes P_n$ & $3 \le md(P_n \boxtimes P_n) \le 4$, $n \ge 7$ \cite{hakanen2024complexity} \\
\hline
Kayak paddle graph $KP(\vartheta, \lambda, \mu)$, $\vartheta, \lambda, \mu \ge 4$ & $3$ \cite{ikhlaq2023new} \\
\hline
Dragon graph $T_{n,m}$, $n \ge 4, m \ge 3$ & $3$ \cite{ikhlaq2023new} \\
\hline
The comb product of $P_n$ and $P_m$, $n,m \ge 4$ & $3$ \cite{ikhlaq2023new} \\
\hline
Tree $T$ of order $n$ and diameter $\ge 2$, $md(T) < \infty$ & $\le n-2$ \cite{hafidh2019multiset} \\
\hline
Prism graph $K_2 \square C_{13}$ & $\ge 4$ \cite{marcelo2025multiset} \\
\hline
Prism graph $C_n \square K_2$, $n \ge 6$ & $=
\begin{cases}
3, & n \ge 11, n \ne 13, \\
4, & n = 8, 9, 10, 13, \\
5, & n = 6, 7.
\end{cases}$ \cite{marcelo2025multiset} \\
\hline
Cylindrical graph $P_m \square C_3$, $n=3, m \ge 3$ & $3$ \cite{marcelo2025multiset} \\
\hline
Cylindrical graph $P_m \square C_n$, $m \ge 6, n \ge 3$ & $\le 4$ \cite{marcelo2025multiset2} \\
\hline
Cylindrical graph $P_m \square C_n$, $m \ge 2, n \ge 8m+1$ & $3$ \cite{marcelo2025multiset2} \\
\hline
\end{tabular}}
\end{table}
\FloatBarrier
\section{Multiset Dimension Variations}

In comparison to the large number of variations to the metric dimension, there is very limited work on the multiset dimension, mainly through four variants, namely the local, outer, local(outer), and edge versions. We have conducted a systemic survey of the known results. 

\subsection{The Outer Multiset Dimension} 

\textit{The outer multiset resolving} set was first introduced by Gil-Pons \textit{et al.}\cite{gil2019distance} in 2019. Instead of requiring that all vertices be distinguished, this modified version requires only vertices \textit{outside} the set $S$ to have distinct multiset representations. Formally, a set $S \subset V(G)$ is an outer multiset resolving set if for any two distinct vertices $u, v \in V(G) \setminus S$, their multisets of distances to $S$ differ \cite{klavvzar2023further}. This approach guaranties the existence of such sets for every graph and allows the definition of the outer dimension of the multiset. This makes it a more practical and computable parameter, especially for graphs where the classical multiset dimension is undefined or infinite. However, further studies on the outer multiset dimension remain limited.

\begin{definition} [Outer Multiset Dimension] \cite{klavvzar2023further} \cite{simanjuntak2025local} 
If for every pair of distinct vertices $u, v \in V(G) \setminus W$, 
$m(u \mid W) \ne m(v \mid W)$, then $W$ is called an \textit{outer multiset resolving set} for $G$. 
This set with the smallest cardinality is an \textit{outer multiset basis}, 
and its cardinality is the \textit{outer multiset dimension} of $G$, 
denoted by $\operatorname{dim}_{ms}(G)$.
\end{definition}

It is obvious that a multiset basis of a graph is also its outer multiset basis. Thus, $\dim(G) \le \dim_{ms}(G) \le md(G)$ \cite{simanjuntak2025local}. Furthremore, we have.

\begin{theorem}\cite{gil2019distance}\cite{klavvzar2023further}\cite{simanjuntak2025local}\label{outertheorem1}
If the order of $G$ is $n$, then $\dim_{ms}(G) \le n - 1$. 
Moreover, $\dim_{ms}(G) = n - 1$ if and only if $G$ is a regular graph 
with diameter at most $2$.    
\end{theorem}

The quantity $f(n,d)$ defined in \cite{chartrand2000resolvability} is the smallest positive integer $k$ such that $k + d^k \ge n$. 
Similarly, $f'(n,d)$ is defined as the smallest positive integer $k'$ such that
$k' + \binom{r+d-1}{d-1} \ge n$. We have the following. 

\begin{theorem}\cite{gil2019distance}
For every graph $G = (V, E)$ of order $n$ and diameter $d$ such that 
$\dim(G) < f'(n,d)$,  
$\dim_{\mathrm{ms}}(G) > \dim(G)$.
\end{theorem}

Two different vertices $u$, $v$ are called true twins if $ N[u] = N[v]$. Likewise, $u$, $v$ are called false twins if $N(u) = N(v)$. In general, $u$, $v$ are called twins if they are either true twins or false twins. 
The property of being twins induces an equivalence relation on the vertex set of any graph \cite{gil2019distance}.

\begin{corollary}\cite{gil2019distance}
Let $G$ be a non-trivial graph and let 
$\mathcal{T} = \{ [u_1], [u_2], \ldots, [u_t] \}$ 
be the set of equivalence classes induced in $V(G)$ by the twin equivalence relation.  
Then, 
$
\dim_{\mathrm{ms}}(G) \ge \sum_{i=1}^{t} (\, |[u_i]| - 1 \,).
$
\end{corollary}

Similarly to the multiset dimension problem, the problem of determining the outer multiset dimension is NP-complete. then the attention is to graphs with small values of $\dim_{\mathrm{ms}}(G)$. We have the following. 

\begin{lemma}\cite{klavvzar2023further}
If $G$ is a graph with $\dim_{\mathrm{ms}}(G) = 2$ and 
$S = \{u, v\}$ is an outer multiset basis, then $d_G(u, v) \le 2$.
\end{lemma}

\begin{theorem}\cite{klavvzar2023further}
Deciding whether a graph $G$ of order $n$ satisfies $\dim_{\mathrm{ms}}(G) = 2$ can be done in $O(n^3)$ time.
\end{theorem}

For a graph $G$, $X \subseteq V(G)$, and $k \ge 0$, define
$
L_k(X) = \{ u \in V(G) : \min_{x \in X} d_G(u, x) = k \}.
$
Note that $L_0(X) = X$, and that the sets $L_k(X)$, $k \ge 0$, partition $V(G)$ \cite{klavvzar2023further}.

\begin{lemma}\cite{klavvzar2023further}
Let $G$ be a graph with $\dim_{\mathrm{ms}}(G) = 2$ and let $S$ be an outer multiset basis of $G$. 
Then for every $k \ge 1$ we have $|L_k(S)| \le 3$. Moreover, if the vertices of $S$ are adjacent, 
then $|L_k(S)| \le 2$ and $|L_{k+1}(S)| \le |L_k(S)|$.
\end{lemma}

Furthermore, in \cite{klavvzar2023further}, the family $\mathcal{F}$ of graphs $G$ is constructed in the following way. 
The vertex set of $G \in \mathcal{F}$ is 
$V(G) = \{u_0, \dots, u_r\} \cup \{v_0, \dots, v_s\}$ for some $r \ge 0$ and $s \ge 1$, 
and the edges of $G$ are given as follows:
\begin{itemize}
    \item $u_0v_0, \ u_0v_1 \in E(G)$.
    \item For every $i \in [r]$ and every $j \in [s]$, $u_{i-1}u_i \in E(G)$ and $v_{j-1}v_j \in E(G)$.
    \item For every $i \in [\min\{r, s\}]$, the edge $u_iv_i$ might exist or not in $G$.
    \item For every $j \in [\min\{r, s-1\}]$, the edge $u_iv_{i+1}$ might exist or not in $G$.
\end{itemize} 

Then the following structural results are shown.
\begin{theorem}\cite{klavvzar2023further}
A graph $G$ has an outer multiset basis formed by two adjacent vertices if and only if $G \in \mathcal{F}$.
\end{theorem}

Recall that the \textit{lexicographic product} $G \circ H$ of graphs $G$ and $H$ has the vertex set $V(G) \times V(H)$ and edges $(g,h)(g',h')$, where either $g = g'$ and $hh' \in E(H)$, or $gg' \in E(G)$. If $g \in V(G)$, then the set of vertices $\{(g,h) : h \in V(H)\}$ induces a subgraph of $G \circ H$ isomorphic to $H$ which is called an $H$-layer and denoted by $^{g}H$. Also, a graph $G$ is\textit{ multiset distance irregular} if for every two vertices $u, v \in V(G)$, the multisets $m_G(u \mid V(G))$ and $m_G(v \mid V(G))$ are different \cite{klavvzar2023further}.

\begin{theorem}\cite{klavvzar2023further}
If $G$ is a graph with $n(G) \ge 2$ and $H \in \{K_k,\overline {K_k}\}$, $k \ge 2$, then $\dim_{ms}(G \circ H) \ge n(G)(k - 1)$ where $n(G)$ is the order of the graph. Moreover, equality holds if and only if $G$ is multiset distance irregular.
\end{theorem}


The following table summarizes the results on some fundamental family of graphs. 

\begin{table}[!ht]
\centering
\caption{The outer multiset dimension of some graphs}
\label{tab:Omd}
\renewcommand{\arraystretch}{1.25}
\setlength{\arrayrulewidth}{0.5pt}
\resizebox{\textwidth}{!}{
\begin{tabular}{|l|c|}
\hline
\textbf{Graph} & \boldmath$\dim_{ms}$  \\
\hline

Path graph $P_n$ &$1$ \cite{gil2019distance}\\
\hline

Cycle graph $C_4$ & $3$ \cite{gil2019distance} \\
\hline

Cycle graph $C_5$ & $4$ \cite{gil2019distance} \\
\hline

Cycle graph $C_n$ & $3$, $n\ge6$ \cite{gil2019distance} \\
\hline

Complete graph $K_n$, $n \ge 2$ & $n - 1$ \cite{alfarisi2019local} \\
\hline

Complete $k$-partite ($K_{r_1,r_2,\dots,r_k}$,  $r_1 = r_2 = \dots = r_k \ge 2$) & $n - 1$ \cite{gil2019distance} \\
\hline

$(P_s \,\square\, P_t)$, $s \ge t \ge 2$  & $3$\cite{gil2019distance}\cite{klavvzar2023further} \\
\hline

\end{tabular}}
\end{table}
\FloatBarrier


\subsection{Local Multiset Dimension} 
\noindent
\par

In \cite{alfarisi2019local}, the authors introduced  the concept of \textit{local multiset dimension} 
$\mu_{\ell}(G)$ (also written as ${lmd}(G)$), which focuses on distinguishing only adjacent vertices using multiset representations. They studied this parameter on several fundamental families of graphs, including paths, cycles, complete graphs, and trees. Their results provide exact values and sharp bounds that establish a baseline understanding of 
${lmd}(G)$. The definition of the local multiset dimension is given below.

\begin{definition} [Local Multiset Dimension] \cite{alfarisi2019local}
Let $G$ be a connected graph with vertex set $V(G)$. Suppose $W = \{s_1, s_2, \ldots, s_k\}$ of vertex set $G$, the multiset representation
of a vertex $v$ of $G$ with respect to $W$ is $r_m(v|W) = \{d(v, s_1), d(v, s_2), \ldots, d(v, s_k)\}$ where $d(v, s_i)$ is a distance between of the vertex $v$
and the vertices in $W$ together with their multiplicities. The resolving set $W$ is a \textit{local resolving set} of $G$ if $r_m(v|W) \neq r_m(u|W)$ for every pair $u, v$ of adjacent vertices of $G$. The minimum local resolving set $W$ is a \textit{local multiset basis} of $G$. If $G$ has a local multiset basis,
then its cardinality is called \textit{local multiset dimension}, denoted by $\mu_{l}(G)$, ${lmd}(G)$ or ${md}_l(G)$.
\end{definition}

Every local multiset basis of a graph induces its local basis, 
and every multiset basis of a graph induces its local multiset basis. 
This leads to the relationship ${ldim}(G)$ $\le$${lmd}(G)$ $ \le$ ${md}(G)$ \cite{alfarisi2019local}.

In a later paper \cite{alfarisi2020note}, the authors discussed the relationship between the multiset dimension and local multiset dimension of graphs. The key insights is that the multiset dimension and local multiset dimension are not monotonic with respect to the number of vertices or the number of edges of a graph, also the difference between the two parameters can be arbitrarily large i.e., gap$(md(G), {lmd}(G)) = \infty$. The investigation of the local multiset dimension was also extended to graphs obtained via a standard graph operation, in particular the Cartesian product and several bounds was established.

\begin{lemma} \cite{alfarisi2020note}
Given that two connected graphs $G_1$ and $G_2$, ${lmd}(G_1 \square G_2) \geq \min \{{lmd}(G_1), {lmd}(G_2)\}$.  
\end{lemma}

More specific, the Cartesian product of a graph $G$ and a tree graph $T$, where ${lmd}(G) = 1$, yields the following results:

\begin{theorem} \cite{alfarisi2020note}
Given that a connected graph $G$ and a path $P_n$, ${lmd}(G \square P_n) = {lmd}(G)$.  
\end{theorem}

\begin{theorem} \cite{alfarisi2020note}
Given that a connected graph $G$ and a tree $T$, ${lmd}(G \square T) \leq {lmd}(G)$.   
\end{theorem}

\begin{corollary} \cite{alfarisi2020note}
Given that a connected graph $G$ and a tree $T$. For ${lmd}(G) = 1$, ${lmd}(G \square T) = 1$. 
\end{corollary}

\begin{theorem} \cite{alfarisi2020note}
For ${lmd}(G) = 1$ and $m$ is even, ${lmd}(G \square C_m) = 1$. 
\end{theorem}

\begin{lemma} \cite{alfarisi2020note}
For $({lmd}(G) = 1 \text{ and } m \text{ is odd}$ or ${lmd}(G) \neq 1 \text{ and } m \geq 3)$, ${lmd}(G \square C_m) \geq 2$.    
\end{lemma}

And in the case where ${lmd}(G) \neq 1$, we have the following. 

\begin{lemma} \cite{alfarisi2020note}
For ${lmd}(G) \neq 1$ and any tree $T$, ${lmd}(G \square T) \geq 2$.   
\end{lemma}

In a recent study \cite{alfarisi2025local}, the comb product of two connected graphs, demoted by $\triangleright$, was introduced. In the following, the main results on the local multiset dimension of comb product graphs with respect to the main graph is presented.

\begin{lemma} \cite{alfarisi2025local}
Let $G$ and $H$ be a nontrivial connected graph.
$lmd(G) = lmd(H) = 1$ \text{if and only if} $ lmd(G \triangleright H) = 1.$
\end{lemma}

\begin{theorem} \cite{alfarisi2025local}
Let $T$ be a tree with order $n \ge 2$, and let $T_1, T_2$ be subtrees of $T$ with
$V(T_1) \cup V(T_2) = V(T)$ and $V(T_1) \cap V(T_2) = \emptyset$, then
\[
lmd(T \triangleright C_3) \ge
\begin{cases}
n, & \text{for } T_1 \ne T_2,\\
n+1, & \text{for } T_1 = T_2.
\end{cases}
\]
\end{theorem}

The Table \ref{tab:lmd} summarize all the known results on some families of graphs. The \textit{$m$-shadow} of a connected graph $G$, denoted by $D_m(G)$, is constructed by taking $m$ copies of $G$, say $G_1, G_2, \odot, G_m$, then joining each vertex $u$ in $G_i$ to the neighbors of the corresponding vertex $v$ in $G_j$, where $1 \le i$ and $j \le m$ \cite{adawiyah2019localmshadow}.

\begin{definition}\cite{alfarisi2024local} 
Let $(G_i)$ be a finite collection of graphs and each $G_i$ has a fixed vertex $v$ called a terminal. The amalgamation ${Amal}(G_i, v, m)$ is formed by taking all the $G_i$s and joining them on the terminal.
\end{definition}

\begin{lemma}\cite{alfarisi2024local}
Let $m$ and $n$ be two integers with $m \ge 2$ and $n \ge 3$. Let $G$ be a connected graph of order $n$ and $v$ be a terminal vertex of $G$, then 
$
{lmd}({Amal}(G, v, m)) \le m \cdot {lmd}(G).
$
\end{lemma}

 \begin{table*}[!ht]
\centering
\caption{Local multiset dimension of some graphs}
\label{tab:lmd}
\renewcommand{\arraystretch}{1.4}
\setlength{\arrayrulewidth}{0.4pt}
\resizebox{\textwidth}{!}{
\begin{tabular}{|l|c|}
\hline
\textbf{Graph} & \boldmath$lmd$  \\
\hline

Path graph $P_n$ &$1$, $n \ge 3$ \cite{alfarisi2019local}\\
\hline

Path graph with homogeneous pendant edges $P_n \,\odot\, mK_1$ &$1$, $n \ge 3$ \cite{adawiyah2020}\\
\hline

$P_n \triangleright C_m$, $n \ge 2$ and $m \ge 3$ & 
$\displaystyle \begin{cases}
1, & \text{for } m \text{ is even},\\
n, & \text{for } m \text{ is odd and } m \ne 3 \text{ or for } m=3 \text{ and } n \text{ is odd},\\
n+1, & \text{for } m = 3 \text{ and } n \text{ is even}.
\end{cases}$ \cite{alfarisi2023local}\\
\hline

Tree $T$ of order $n$ & $1$ \cite{alfarisi2024local}\cite{alfarisilocal}\\
\hline

Complete $k$-ary tree of height $h$& $1$ \cite{alfarisi2019local}\\
\hline

$T_1 \triangleright T_2$ & $1$  \cite{alfarisi2023local}\\
\hline

$T \triangleright C_m$, $m \ge 4$ & $1$ for $m$ even \cite{alfarisi2023local}\\
\hline

Caterpillar graph $C_{n,m}$ & $1$,  $n \ge 3$ and $m \ge 1$ \cite{alfarisi2019local}\\
\hline

Caterpillar with adding one edge $C'_{n,m}$& $\infty$, $n \ge 3$ \cite{adawiyah2019}\\
\hline

Caterpillar graph with homogeneous pendant edges $C_{n,m} \,\odot\, mK_1$ & $1$, $n \ge 3$ and $m \ge 2$ \cite{adawiyah2020}\\
\hline

Bipartite graph & $1$ \cite{alfarisi2024local}\cite{alfarisilocal}\\
\hline

$k$-partite graph $K_{n_1, n_2, \ldots, n_k}$  with $1 \le l \le k$ and $n_l \ge k - 1$ & $\frac{k(k - 1)}{2}$ \cite{alfarisi2019local}\\
\hline

Star graph $S_n$ & $1$, $n \ge 3$ \cite{alfarisi2019local}\\
\hline

Star graph with homogeneous pendant edges $S_n \,\odot\, mK_1$ & $1$, $n \ge 3$ \cite{adawiyah2020}\\
\hline

Cycle graph $C_n$ with $n \ge 3$ & 
$\begin{cases}
1, & \text{if } n \text{ is even},\\[4pt]
3, & \text{if } n \ge 7 \text{ is odd},\\
\infty, & \text{if } n=3,5.
\end{cases}$
\cite{alfarisi2019local}\cite{simanjuntak2025local}\\
\hline

Cycle graph with two pendant $C^2_n$, $n \ge 4$, $4 \le i \le n$ & 
$\begin{cases}
1, & \text{if } m \equiv 0 \pmod{2}\\
2, & \text{if } m \equiv 1 \pmod{2}.
\end{cases}$ 
\cite{adawiyah2019}\\
\hline

Cycle graph with homogeneous pendant edges $C_n \,\odot\, mK_1$, $n \ge 4$ & 
$\begin{cases}
1, & \text{if } n \text{ is even},\\
2, & \text{if } n \text{ is odd}.
\end{cases}$
\cite{adawiyah2020}\\
\hline

$C_n \triangleright T$, $n \ge 3$ & $\begin{cases}
1, & n \text{ even},\\
2, & n \text{ odd}.
\end{cases}$ \cite{alfarisi2023local}\\
\hline

$C_n \triangleright C_m $, $n,m \ge 3$ & $\begin{cases}
1, & \text{both } n \text{ and } m \text{ are even},\\
2, & \text{one of } n \text{ and } m \text{ is even and the other is odd},\\
n, & \text{both } n \text{ and } m \text{ are odd}.
\end{cases}$ \cite{alfarisi2023local}\\
\hline

Complete graph $K_2$& 
$1$ \cite{alfarisi2019local}\\
\hline

Complete graph $K_n$& 
$\infty$, $n \ge 3$
\cite{alfarisi2019local}\\
\hline

Complete graph with homogeneous pendant edges $K_n \,\odot\, mK_1$& 
$\infty$, $n \ge 3$
\cite{adawiyah2020}\\
\hline

Wheel graph $W_n$, $n \ge 3$ &
$\begin{cases}
3, & n = 4,6,\\
\left\lceil n/4 \right\rceil, & n \ge 8 \text{ and even},\\
\infty, & \text{otherwise.}
\end{cases}$
\cite{simanjuntak2025local}\\
\hline

Pan graph $P'_n$, $n \ge 4$, $4 \le i \le n$ & 
$\begin{cases}
1, & \text{if } n \equiv 0 \pmod{2}\\
2, & \text{if } n \equiv 1 \pmod{2}.
\end{cases}$
\cite{adawiyah2019}\\
\hline

Sunlet graph $Sl_{_n}$, $n \ge 4$, $4 \le i \le n$& 
$\begin{cases}
1, & \text{if } m \equiv 0 \pmod{2}\\
2, & \text{if } m \equiv 1 \pmod{2}.
\end{cases}$
\cite{adawiyah2019}\\
\hline

Non-cycle unicyclic graph & $\begin{cases}
1, & \text{if } G \text{ contains an even cycle},\\
2, & \text{if } G \text{ contains an odd cycle.}
\end{cases}$ \cite{simanjuntak2025local}\\
\hline

Unicyclic graph of order at least 3 and contains a circle of order $p$& 
$\begin{cases}
1, & p \text{ even},\\
2, & p \text{ odd.}
\end{cases}$ \cite{alfarisi2023local}\\ 
\hline

Bicyclic graph of order at least 3 and contains cycles of order $p_1$ and $p_2$ & 
$\begin{cases}
1, & \text{if } p_1 \text{ and } p_2 \text{ are even},\\
2, & \text{otherwise}.
\end{cases}$ \cite{alfarisi2023local}\\
\hline

m-shadowing path graph $D_m(P_n)$ & $1$, $n \ge 3$ \cite{adawiyah2019localmshadow}\\
\hline

m-shadowing cycle graph $D_m(C_n)$ & $\begin{cases}
1, & \text{if } n \text{ is even},\\
3, & \text{if } n \text{ is odd}.
\end{cases}, n \ge 3$ \cite{adawiyah2019localmshadow}\\
\hline

m-shadowing star graph $D_m(S_n)$ & $1$, $n \ge 3$ \cite{adawiyah2019localmshadow}\\
\hline

m-shadowing gear graph $D_m(G_n)$ & $1$, $n \ge 4$ \cite{adawiyah2019localmshadow}\\
\hline

m-shadowing complete graph $D_m(K_n)$ & $\infty$, $n \ge 3$ \cite{adawiyah2019localmshadow}\\
\hline

Amalgamation of $m$ paths ${Amal}(P_n, v, m)$ & $1$ \cite{alfarisi2024local}\\
\hline

Amalgamation of $m$ complete graphs ${Amal}(K_n, v, m)$ & $\infty$ \cite{alfarisi2024local}\\
\hline

Amalgamation of $m$ wheel graph ${Amal}(W_n, v, m)$ ($m \ge 2$, $n \ge 3$, $v$ center) & 
$ m \cdot {lmd}(W_n)$ \cite{alfarisi2024local} \\
\hline

Amalgamation of $m$ wheel graph ${Amal}(W_n, v, m)$  ($v$ not center) &
$\begin{cases}
m\cdot\left(\frac{n}{4}-1\right), & n \equiv 0 \pmod{4},\\[0.3em]
m\cdot\left\lfloor n/4 \right\rfloor, & n \equiv 1,2,3 \pmod{4}.
\end{cases}$
\cite{alfarisi2024local} \\
\hline

Amalgamation of $m$ fan graph ${Amal}(F_n, v, m)$ ($m \ge 2$, $n \ge 6$, $v$ center) & 
$ m \cdot {lmd}(F_n)$ \cite{alfarisi2024local} \\
\hline

Amalgamation of $m$ fan graph ${Amal}(F_n, v, m)$ ($m \ge 2$, $n \ge 3$, $d(v)=2$) & 
$m \lfloor n/4 \rfloor$ \cite{alfarisi2024local} \\
\hline

Amalgamation of $m$ fan graph ${Amal}(F_n, v, m)$ ($m \ge 2$, $n \ge 3$, $d(v)=3$) &
$\displaystyle m\left( \Big\lfloor \frac{p}{4} \Big\rfloor + \Big\lfloor \frac{n-p}{4} \Big\rfloor \right) \text{ or } m\left( \Big\lfloor \frac{p}{4} \Big\rfloor + \Big\lceil \frac{n-p}{4} \Big\rceil \right)$ \cite{alfarisi2024local}\\
\hline

\end{tabular}}
\end{table*}
\FloatBarrier
\subsection{Local Outer Multiset Dimension}
The local version of the outer multiset dimension was introduced in \cite{simanjuntak2025local}.

\begin{definition}\cite{simanjuntak2025local} 
 If for every pair of distinct vertices $u, v \in V(G) \setminus W$, 
$m(u \mid W) \ne m(v \mid W)$, then $W$ is called a \textit{local outer multiset resolving set} for $G$. 
This set with the smallest cardinality is a \textit{local outer multiset basis}, 
and its cardinality is the \textit{local outer multiset dimension} of $G$, 
denoted by ${ldim}_{ms}(G)$.   
\end{definition}

Since every outer multiset basis of a graph induces its local outer multiset basis, 
then ${ldim}_{ms}(G) \le \dim_{ms}(G)$. 
Additionally, since every local multiset basis of a graph induces its local outer multiset basis, 
and every local outer multiset basis of a graph induces its local basis, 
then ${ldim}(G) \le {ldim}_{ms}(G) \le {lmd}(G)$\cite{simanjuntak2025local}.

The paper \cite{simanjuntak2025local} established fundamental results on the local outer multiset dimension, including key properties, parameter relations, structural bounds, and results for several graph families such as amalgamations of complete graphs, non-cycle unicyclic graphs, cycles, and wheel graphs.

\begin{table}[!ht]
\centering
\caption{Local outer multiset dimension of some graphs}
\label{tab:LOmd}
\renewcommand{\arraystretch}{1.25}
\setlength{\arrayrulewidth}{0.4pt} 
\resizebox{\textwidth}{!}{
\begin{tabular}{|l|c|}
\hline
\textbf{Graph} & $\boldsymbol{ldim_{ms}}$  \\
\hline

Amalgamation of complete graphs $\text{Amal}(K_{n_i}, m)$ 
& $\begin{cases}
1, & m_3 = 0,\; m_{\ge 4} = 0,\\
2, & m_3 = 1,\; m_{\ge 4} = 0,\\
\displaystyle \sum_{n_i \ge 3} (n_i - 2), & \text{otherwise.}
\end{cases}$
\cite{simanjuntak2025local}\\ 
\hline

Bipartite graph & $1$ \cite{simanjuntak2025local}\\
\hline

Non-cycle unicyclic graph & $\begin{cases}
1, & \text{if $G$ contains an even cycle},\\
2, & \text{if $G$ contains an odd cycle.}
\end{cases}$ \cite{simanjuntak2025local}\\
\hline

Cycle graph $C_n$, $n \ge 3$ &
$\begin{cases}
1, & \text{if $n$ is even},\\
2, & \text{if $n$ is odd}.
\end{cases}$
\cite{simanjuntak2025local}\\
\hline

Complete graph $K_n$, $n \ge 2$ & $n - 1$ \cite{alfarisi2019local} \\
\hline

Wheel graph $W_n$, $n \ge 3$ &
$\begin{cases}
3, & n = 3,4,6,\\
\left\lceil n/4 \right\rceil, & n \ge 8 \text{ even or } n \equiv 1 \pmod{4},\\
\left\lceil n/4 \right\rceil + 1, & \text{otherwise.}
\end{cases}$ \cite{simanjuntak2025local}\\
\hline

\end{tabular}}
\end{table}
\FloatBarrier


\subsection{Edge Multiset Dimension}

\noindent
\par

\vspace{10pt}
The concept of the edge multiset dimension of graphs was introduced in \cite{ikhlaq2023new}.

\begin{definition}\cite{ikhlaq2023new} 
Let $G$ be a graph and let $B = \{b_1, b_2, \ldots, b_j\}$ be a set of vertices. 
For an edge $e = xy$, the \emph{edge-multiset representation} of $e$ with respect to $B$ 
is the multiset
\[
{r_{em}}(e \mid B) = 
\{\, d_{G}(e, b_1),\, d_{G}(e, b_2),\, \ldots,\, d_{G}(e, b_j) \,\}.
\]
If for any two distinct edges $e$ and $f$ of $G$ we have 
${r_{em}}(e \mid B) \neq {r_{em}}(f \mid B)$, 
then $B$ is called an \emph{edge-multiset resolving set} for $G$.
If $G$ admits an edge-multiset resolving set, then a smallest such set is called 
an \emph{edge-multiset basis} of $G$.  
The number of elements in an edge-multiset basis is called the 
edge-multiset dimension of $G$, denoted by ${md}_e(G)$.
If $G$ does not contain any edge-multiset resolving set, we define 
${md}_e(G) = \infty$.
\end{definition}

In \cite{ikhlaq2023new}, it was shown that the edge multiset dimension of a connected graph $G$ is at least the edge metric
dimension of graph $G$, that is, $\operatorname{edim}(G) \le md_e(G)$ . Moreover, no connected graph has edge multiset dimension equal to $2$ ,i.e., $md_e(G) \neq 2$. 
In particular, for any connected graph other than a path, $md_e(G) \ge 3$. Several exact values of this parameter for basic graph classes were also determined, see Table~ \ref{tab:emd}.

Several structural conditions have been identified that ensure a graph has infinite edge multiset dimension:
\begin{lemma}\cite{ikhlaq2023new}
Let $G$ be a graph with a set of vertices and edges, $V$ and $E$, respectively, where $|V| \ge 2$.
If the distance between a vertex and an edge is at most $2$, then $G$ does not contain an edge-multiset
resolving set.
\end{lemma}

\begin{theorem}\cite{ikhlaq2023new}
If any graph $G$ has a vertex adjacent to at least three pendant vertices, then $md_e(G) = \infty$.    
\end{theorem}

The following families of graphs have infinite edge multiset dimensions. Complete
graph, star graph, friendship graph, wheel graph, the Peterson graph, fan graph, complete bipartite
graph and cycle graph with at most 6 vertices \cite{ikhlaq2023new}.

Some bounds on the edge multiset dimension of trees:
\begin{lemma}\cite{ikhlaq2023new}
Let $T_n$ be a tree with order $n$ and diameter $3$. If $md_e(T_n) \neq \infty$, then $md_e(T_n) \le n - 2$.
\end{lemma}

\begin{theorem}\cite{ikhlaq2023new}
Let $T_n$ be a tree with order $n$ and $T_n$ as not a path. If $md_e(T_n) \neq \infty$, then $3 \le
md_e(T_n) \le n - 2$.
\end{theorem}


\noindent
\par



\begin{table}[!ht]
\centering
\caption{Edge-multiset dimension of some graphs}
\label{tab:emd}
\renewcommand{\arraystretch}{1.25}
\setlength{\arrayrulewidth}{0.8pt} 
\resizebox{\textwidth}{!}{ \large
\begin{tabular}{|l|c|}
\hline
\textbf{Graph} & $\mathbf{md_e}$ \\
\hline

Path graph $P_n$ & $1$\cite{ikhlaq2023new}\\
\hline

Complete graph & $\infty$\cite{ikhlaq2023new}\\
\hline

Star graph & $\infty$\cite{ikhlaq2023new}\\
\hline

Friendship graph & $\infty$\cite{ikhlaq2023new}\\
\hline

Wheel graph & $\infty$\cite{ikhlaq2023new}\\
\hline

The Peterson graph & $\infty$\cite{ikhlaq2023new}\\
\hline

Fan graph & $\infty$\cite{ikhlaq2023new}\\
\hline

Complete bipartite graph & $\infty$\cite{ikhlaq2023new}\\
\hline

Cycle graph with at most 6
vertices & $\infty$\cite{ikhlaq2023new}\\
\hline

Cycle graph $C_n$, $n \ge 7$ & $3$ \cite{ikhlaq2023new}\\
\hline

Kayak paddle graph $KP(\vartheta, \lambda, \mu)$, \text{with } $\vartheta, \lambda, \mu \ge 4$ & $3$ \cite{ikhlaq2023new}\\
\hline

Dragon graph $T_{n,m}$ with $n \ge 4$ and $m \ge 3$ & $3$ \cite{ikhlaq2023new}\\
\hline

The comb product of $P_n$ and $P_m$ with $n, m \ge 4$ & $3$ \cite{ikhlaq2023new}\\
\hline

Caterpillar graph $CT_n$,$n \ge 5$ & $3$ 
\cite{ikhlaq2023new}\\
\hline

Lobster graph $L_n$, $n \ge 5$ & $3$ \cite{ikhlaq2023new}\\
\hline

Caterpillar graph $CT_n$ with $n \ge 5$, $\deg(x_i) \le 4$ where $3 \le i \le n - 2$ and $D = \{i \mid 3 \le i \le n - 2\}$&  $3 + |D|$ \cite{ikhlaq2023new}\\
\hline

Lobster graph $L_n$ with $n > 5$, $\deg(x_i) = 4$, where $i \in D = \{i \mid 3 \le i \le n - 2\}$ & $ 3 +|D|$ \cite{ikhlaq2023new}\\
\hline

\end{tabular}}
\end{table}
\FloatBarrier

\section{New Directions in Multiset-Based Metric Dimensions}

From the above survey, we see many variants have been developed based on the classical metric dimension, but their multiset version remain limited, thus, we recommend investigating some new variations to the multiset metric dimension.

\subsection{Multiset Partition Dimension}

Based on the partition dimension of a graph \cite{chartrand2000partition}, we introduce the following parameter.   

Let $G = (V,E)$ be a connected graph and let 
$\Pi = \{P_1, P_2, \dots, P_k\}$ be an ordered partition of $V(G)$. 
For a vertex $v \in V(G)$, the \emph{multiset representation} of $v$ with respect to $\Pi$ is defined as
\[
r_m(v \mid \Pi) = \{\, d(v, P_1), d(v, P_2), \dots, d(v, P_k) \,\},
\]
where $d(v, P_i) = \min\{ d(v,x) : x \in P_i \}$. 
The values are taken as a multiset, allowing repetition but ignoring order.

The partition $\Pi$ is called a \emph{multiset resolving partition} if for every pair of distinct vertices $u, v \in V(G)$,
\[
r_m(u \mid \Pi) \neq r_m(v \mid \Pi).
\]
The \emph{multiset partition dimension} of $G$, denoted by $\operatorname{mpd}(G)$, 
is the minimum cardinality of such a partition.

Apparently, this parameter behaves quite differently. In \cite{bong2021some}, a graph with $md=3$ with the largest possible number of vertices have been given, however, the $mpd$ here is much larger.

\begin{figure}[h!]
    \centering
    \includegraphics[width=0.50\textwidth]{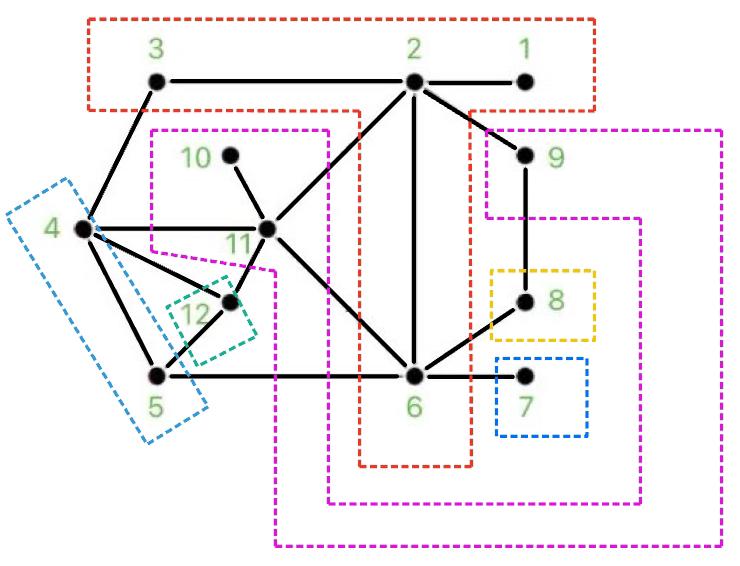}
    \caption{The graph has $\operatorname{mpd}=6$, with partitions illustrated in the figure.}
    \label{fig:Novigraphwithmpd=6}
\end{figure}
\FloatBarrier

\begin{table}[!ht]
\centering
\scriptsize
\begin{tabular}{ll}
\toprule
\textbf{Vertex} & 
\textbf{Multiset partition representation}\\
\midrule
1 & $\{0, 2, 3, 3, 3, 3\}$ \\
2 & $\{0, 1, 2, 2, 2, 2\}$ \\
3 & $\{0, 1, 2, 2, 3, 3\}$ \\
4 & $\{0, 1, 1, 1, 3, 3\}$ \\
5 & $\{0, 1, 1, 2, 2, 2\}$ \\
6 & $\{0, 1, 1, 1, 1, 2\}$ \\
7 & $\{0, 1, 2, 2, 2, 3\}$ \\
8 & $\{0, 1, 1, 2, 2, 3\}$ \\
9 & $\{0, 1, 1, 3, 3, 3\}$ \\
10 & $\{0, 2, 2, 2, 3, 3\}$ \\
11 & $\{0, 1, 1, 1, 2, 2\}$ \\
12 & $\{0, 1, 1, 2, 3, 3\}$ \\
\bottomrule
\end{tabular}
\caption{Multiset partition representations of the vertices with respect to the partition shown in Figure~\ref{fig:Novigraphwithmpd=6}.}
\end{table}

It is natural to investigate how the multiset partition dimension can be bounded in terms of the multiset dimension and other metric type parameters. In particular, establishing tight upper and lower bounds is essential for understanding its behavior and for identifying extremal graph classes. It is also interesting to examine whether the multiset partition dimension is a characteristic parameter, in the sense that it captures intrinsic structural properties of graphs or distinguishes between non-isomorphic graphs in a meaningful way. 



\subsection{Other Multiset-based Variants}

Analogously to the concept of variations to the metric dimension of graphs, we introduce the following variations.

\begin{definition} 
Let $G$ be a graph and let $B = \{b_1, b_2, \ldots, b_j\}$ be a set of vertices. The Edge Multiset Representation of $e$ with respect to $B$ 
is the multiset
\[
{r_{em}}(e \mid B) = 
\{\, d_{G}(e, b_1),\, d_{G}(e, b_2),\, \ldots,\, d_{G}(e, b_j) \,\}.
\]
If any two incident edges, i.e. the edges sharing the common end vertex, have different representation, then $B$ is called an \emph{local edge-multiset resolving set} for $G$. The smallest such set is called an \emph{local edge multiset basis} of $G$.  
The number of elements in an edge-multiset basis is called the 
local edge multiset dimension of $G$, denoted by ${md}_e(G)$.
\end{definition}

\begin{definition} 
Let $G=(V,E)$ be a simple connected graph and let $S \subseteq V$. For each vertex $v \in V$, the multiset representation of $v$ with respect to $S$ is defined as  
\[
r_m(v \mid S) = \{\, d(v,w) : w \in S \,\}.
\]
For any pair of distinct vertices $u,v\in V$, if the representation difference, i.e. the set difference, exceeds k, the set $S$ is called a \textit{k-multiset resolving set} for $G$ and its cardinality is called the \textit{k-multiset dimension} of $G$, denoted by $md_k(G)$.
\end{definition}

It is also possible to consider other multiset-based variants, for example the  \textit{k-multiset antidimension}. However, in the absence of clear theoretical significance or practical applications, we do not formally introduce these parameters in this work.







\begin{thebibliography}{00}

\bibitem{slater1975leaves}
P. J. Slater,
``Leaves of trees,''
\emph{Congr. Numer.}, vol. 14, pp. 549--559, 1975.

\bibitem{harary1976metric}
F. Harary and R. A. Melter,
``On the metric dimension of a graph,''
\emph{Ars Combin.}, vol. 2, pp. 191--195, 1976.

\bibitem{erdos1963two}
P. Erdos and A. R{\'e}nyi,
``On two problems of information theory,''
\emph{Magyar Tud. Akad. Mat. Kutat{\'o} Int. K{\"o}zl}, vol. 8, no. 1--2, pp. 229--243, 1963.

\bibitem{chartrand2000resolvability}
Chartrand, Gary, Eroh, Linda, Johnson, Mark A, and Oellermann, Ortrud R,
``Resolvability in graphs and the metric dimension of a graph,''
\emph{Discrete Applied Mathematics}, vol. 105, no. 1--3, pp. 99--113, 2000, Elsevier.

\bibitem{ikhlaq2023new}
Ikhlaq, Hafiz Muhammad and Ismail, Rashad and Siddiqui, Hafiz Muhammad Afzal and Nadeem, Muhammad Faisal,
``A new technique to uniquely identify the edges of a graph,''
\emph{Symmetry}, vol. 15, no. 3, pp. 762, 2023, MDPI.

\bibitem{bondy2008graph}
Bondy, J. A. and Murty, U. S. R.,
\emph{Graph Theory},
Springer, New York, 2008, ISBN: 978-1-84628-969-9.

\bibitem{kuziak2021metric}
Kuziak, Dorota and Yero, Ismael G,
``Metric dimension related parameters in graphs: A survey on combinatorial, computational and applied results,''
\emph{arXiv preprint arXiv:2107.04877}, 2021.

\bibitem{tillquist2023getting}
Tillquist, Richard C and Frongillo, Rafael M and Lladser, Manuel E,
``Getting the lay of the land in discrete space: A survey of metric dimension and its applications,''
\emph{SIAM Review}, vol. 65, no. 4, pp. 919--962, 2023, SIAM.

\bibitem{mohamed2023comprehensive}
Mohamed, Basma,
``A comprehensive survey on the metric dimension problem of graphs and its types,''
\emph{International Journal of Theoretical and Applied Mathematics}, vol. 9, no. 1, pp. 1--5, 2023.

\bibitem{kratica2014strong}
Kratica, Jozef and Kova{\v{c}}evi{\'c}-Vuj{\v{c}}i{\'c}, Vera and {\v{C}}angalovi{\'c}, Mirjana and Mladenovi{\'c}, Nenad,
``Strong metric dimension: a survey,''
\emph{Yugoslav Journal of Operations Research}, vol. 24, no. 2, pp. 187--198, 2014.

\bibitem{chartrand2000partition}
Chartrand, Gary and Salehi, Ebrahim and Zhang, Ping,
``The partition dimension of a graph,''
\emph{Aequationes mathematicae}, vol. 59, no. 1, pp. 45--54, 2000, Springer.

\bibitem{kelenc2018uniquely}
Kelenc, Aleksander and Tratnik, Niko and Yero, Ismael G,
``Uniquely identifying the edges of a graph: the edge metric dimension,''
\emph{Discrete Applied Mathematics}, vol. 251, pp. 204--220, 2018, Elsevier.

\bibitem{article}
Nasir, Ruby and Zafar, Sohail and Zahid, Zohaib,
``Edge metric dimension of graphs,''
\emph{Ars Combinatoria -Waterloo then Winnipeg-}, 2018.

\bibitem{bailey2011base}
Bailey, Robert F and Cameron, Peter J,
``Base size, metric dimension and other invariants of groups and graphs,''
\emph{Bulletin of the London Mathematical Society}, vol. 43, no. 2, pp. 209--242, 2011, Wiley Online Library.

\bibitem{feng2013metric}
Feng, Min and Xu, Min and Wang, Kaishun,
``On the metric dimension of line graphs,''
\emph{Discrete Applied Mathematics}, vol. 161, no. 6, pp. 802--805, 2013, Elsevier.

\bibitem{caceres2007metric}
C{\'a}ceres, Jos{\'e} and Hernando, Carmen and Mora, Merce and Pelayo, Ignacio M and Puertas, Mar{\'\i}a L and Seara, Carlos and Wood, David R,
``On the metric dimension of cartesian products of graphs,''
\emph{SIAM journal on discrete mathematics}, vol. 21, no. 2, pp. 423--441, 2007, SIAM.

\bibitem{eroh2012metric}
Eroh, Linda and Kang, Cong X and Yi, Eunjeong,
``Metric dimension and zero forcing number of two families of line graphs,''
\emph{arXiv preprint arXiv:1207.6127}, 2012.

\bibitem{okamoto2010local}
Okamoto, Futaba and Phinezy, Bryan and Zhang, Ping,
``The local metric dimension of a graph,''
\emph{Mathematica Bohemica}, vol. 135, no. 3, pp. 239--255, 2010, Institute of Mathematics, Academy of Sciences of the Czech Republic.

\bibitem{fancy2021local}
Fancy, VF and Cynthia, V Jude Annie,
``Local metric dimension of certain wheel related graphs,''
\emph{Int. J. Math. Comput. Sci}, vol. 16, no. 4, pp. 1303--1315, 2021.

\bibitem{alfarisi2019local}
Ridho Alfarisi and Dafik and Arika Indah Kristiana and Ika Hesti Agustin,
``The local multiset dimension of graphs,''
\emph{International Journal of Engineering \& Technology}, vol. 8, no. 3, pp. 120--124, 2019, Science Publishing Corporation, URL: http://www.sciencepubco.com/index.php/IJET.

\bibitem{kelenc2017mixed}
Kelenc, Aleksander and Kuziak, Dorota and Taranenko, Andrej and Yero, Ismael G,
``Mixed metric dimension of graphs,''
\emph{Applied Mathematics and Computation}, vol. 314, pp. 429--438, 2017, Elsevier.

\bibitem{jannesari2012metric}
Jannesari, Mohsen and Omoomi, Behnaz,
``The metric dimension of the lexicographic product of graphs,''
\emph{Discrete mathematics}, vol. 312, no. 22, pp. 3349--3356, 2012, Elsevier.

\bibitem{sebHo2004metric}
Seb{\H{o}}, Andr{\'a}s and Tannier, Eric,
``On metric generators of graphs,''
\emph{Mathematics of Operations Research}, vol. 29, no. 2, pp. 383--393, 2004, INFORMS.

\bibitem{yero2013strong}
Yero, Ismael Gonz{\'a}lez,
``On the strong partition dimension of graphs,''
\emph{arXiv preprint arXiv:1312.1987}, 2013.

\bibitem{adar2017k}
Adar, Ron and Epstein, Leah,
``The k-metric dimension,''
\emph{Journal of Combinatorial Optimization}, vol. 34, pp. 1--30, 2017, Springer.


\bibitem{estrada2013k}
Estrada-Moreno, Alejandro and Rodr{\'\i}guez-Vel{\'a}zquez, Juan A and Yero, Ismael G,
``The k-metric dimension of a graph,''
\emph{arXiv preprint arXiv:1312.6840}, 2013.

\bibitem{arumugam2012fractional}
Arumugam, S and Mathew, Varughese,
``The fractional metric dimension of graphs,''
\emph{Discrete Mathematics}, vol. 312, no. 9, pp. 1584--1590, 2012, Elsevier.

\bibitem{kang2013fractional}
Kang, Cong X and Yi, Eunjeong,
``The fractional strong metric dimension of graphs,''
\emph{International Conference on Combinatorial Optimization and Applications}, pp. 84--95, 2013, Springer.

\bibitem{trujillo2016k}
Trujillo-Rasua, Rolando and Yero, Ismael G,
``k-metric antidimension: A privacy measure for social graphs,''
\emph{Information Sciences}, vol. 328, pp. 403--417, 2016, Elsevier.

\bibitem{Rinovia2017multiset}
{Simanjuntak Rinovia} and Siagian, Presli and Vetrik, Tomas,
``The multiset dimension of graphs,''
\emph{arXiv preprint arXiv:1711.00225}, 2017.

\bibitem{Novi2021some}
Bong, Novi H and Lin, Yuqing,
``Some properties of the multiset dimension of graphs.,''
\emph{Electron. J. Graph Theory Appl.}, vol. 9, no. 1, pp. 215--221, 2021.

\bibitem{alfarisi2020note}
Alfarisi, Ridho and Lin, Yuqing and Ryan, Joe and Dafik, Dafik and Agustin, Ika Hesti,
``A note on multiset dimension and local multiset dimension of graphs,''
\emph{Statistics, Optimization \& Information Computing}, vol. 8, no. 4, pp. 890--901, 2020.

\bibitem{hafidh2019multiset}
Hafidh, Yusuf and Kurniawan, Rizki and Saputro, Suhadi and Simanjuntak, Rinovia and Tanujaya, Steven and Uttunggadewa, Saladin,
``Multiset dimensions of trees,''
\emph{arXiv preprint arXiv:1908.05879}, 2019.

\bibitem{marcelo2025multiset}
Marcelo, Reginaldo M and Garciano, Agnes D and Buot, Jude Cabigas and Tolentino, Mark Anthony C,
``Multiset Dimension of Prisms,''
\emph{Communications in Combinatorics and Optimization}, 2025, Azarbaijan Shahid Madani University.

\bibitem{marcelo2025multiset2}
Marcelo, Reginaldo M and Tolentino, Mark Anthony C and Garciano, Agnes D and Buot, Jude C,
``On multiset dimension of cylindrical graphs,''
\emph{J. COMBIN. MATH. COMBIN. COMPUT}, vol. 126, no. 225, pp. 240, 2025.

\bibitem{chartrand2021distance}
Chartrand, Gary and Kono, Yuya and Zhang, Ping,
``Distance vertex identification in graphs,''
\emph{Journal of Interconnection Networks}, vol. 21, no. 01, pp. 2150005, 2021, World Scientific.

\bibitem{kono2021vertex}
Kono, Yuya and Zhang, Ping,
``Vertex identification in trees,''
\emph{Discrete Math. Lett}, vol. 7, no. 66--73, pp. 2022, 2021.

\bibitem{kono2022note}
Kono, Yuya and Zhang, Ping,
``A note on the identification numbers of caterpillars,''
\emph{Discrete Math. Lett}, vol. 8, pp. 10--15, 2022.

\bibitem{kono2022vertex}
Kono, Yuya and Zhang, Ping,
``Vertex identification in grids and prisms,''
\emph{Journal of Interconnection Networks}, vol. 22, no. 02, pp. 2150019, 2022, World Scientific.

\bibitem{marcelo2025vertex}
Marcelo, Reginaldo M and Tolentino, Mark Anthony C and Garciano, Agnes D and Ruiz, Mari-Jo P and Buot, Jude C,
``On the vertex identification spectra of grids,''
\emph{Journal of Interconnection Networks}, vol. 25, no. 01, pp. 2450002, 2025, World Scientific.

\bibitem{cai2025identification}
Cai, Gaixiang and Xiao, Fengru and Yu, Guidong,
``The identification numbers of lollipop graphs,''
\emph{AIMS Mathematics}, vol. 10, no. 4, pp. 7813--7827, 2025, American Institute of Mathematical Sciences.

\bibitem{wang2024graph}
Wang, Runze,
``Graph identification index,''
\emph{arXiv preprint arXiv:2410.07019}, 2024.

\bibitem{hakanen2024complexity}
Hakanen, Anni and Yero, Ismael G,
``Complexity and Equivalency of Multiset Dimension and ID-colorings,''
\emph{Fundamenta Informaticae}, vol. 191, no. 3--4, pp. 315--330, 2024, SAGE Publications Sage UK: London, England.

\bibitem{gil2019distance}
Gil-Pons, Reynaldo and Ram{\'\i}rez-Cruz, Yunior and Trujillo-Rasua, Rolando and Yero, Ismael G,
``Distance-based vertex identification in graphs: The outer multiset dimension,''
\emph{Applied Mathematics and Computation}, vol. 363, pp. 124612, 2019, Elsevier.

\bibitem{klavvzar2023further}
Klav{\v{z}}ar, Sandi and Kuziak, Dorota and Yero, Ismael G,
``Further contributions on the outer multiset dimension of graphs,''
\emph{Results in Mathematics}, vol. 78, no. 2, pp. 50, 2023, Springer.

\bibitem{simanjuntak2025local}
Simanjuntak, Rinovia and Hasan, M Ali and Anggarawan, Muhung,
``Local (Outer) Multiset Dimensions of Graphs,''
\emph{arXiv preprint arXiv:2507.15071}, 2025.

\bibitem{alfarisi2025local}
Alfarisi, Ridho and Susilowati, Liliek and Kristiana, Arika Indah,
``On the Local Multiset Dimension of Comb Product Graphs,''
\emph{Statistics, Optimization \& Information Computing}, vol. 14, no. 3, pp. 1356--1361, 2025.

\bibitem{adawiyah2019localmshadow}
Adawiyah, R and Agustin, IH and Prihandini, RM and Alfarisi, R and Albirri, ER and others,
``On the local multiset dimension of m-shadow graph,''
\emph{Journal of Physics: Conference Series}, vol. 1211, no. 1, pp. 012006, 2019, IOP Publishing.

\bibitem{alfarisi2024local}
Alfarisi, Ridho and Susilowati, Liliek and Dafik, Dafik and Prabhu, Savari,
``Local multiset dimension of amalgamation graphs,''
\emph{F1000Research}, vol. 12, pp. 95, 2024.

\bibitem{adawiyah2020}
Adawiyah, R and Agustin, IH and Prihandini, RM and Alfarisi, R and Albirri, ER and others,
``On the local multiset dimension of graph with homogenous pendant edges,''
\emph{Journal of Physics: Conference Series}, vol. 1538, no. 1, pp. 012023, 2020, IOP Publishing.

\bibitem{alfarisi2023local}
Alfarisi, Ridho and Susilowati, Liliek and DAFIK, OSAYE J and Osaye, FJ,
``On the local multiset dimension of some families of graphs,''
\emph{WSEAS Trans. Math}, vol. 22, pp. 64--69, 2023.

\bibitem{alfarisilocal}
Alfarisi, R and SusilowatiDafik, L,
``The Local Multiset Resolving of Graphs. 2022,''
Review.

\bibitem{adawiyah2019}
Adawiyah, Robiatul and Prihandini, RM and Albirri, ER and Agustin, IH and Alfarisi, R and others,
``The local multiset dimension of unicyclic graph,''
\emph{IOP Conference Series: Earth and Environmental Science}, vol. 243, no. 1, pp. 012075, 2019, IOP Publishing.

\bibitem{bong2021some}
Bong, Novi H and Lin, Yuqing,
``Some properties of the multiset dimension of graphs.,''
\emph{Electron. J. Graph Theory Appl.}, vol. 9, no. 1, pp. 215--221, 2021.

\end{thebibliography}
\end{document}